\documentclass[10pt,pre,twocolumn,longbibliography,superscriptaddress]{revtex4}
\usepackage[T1]{fontenc}
\usepackage[utf8]{inputenc}
\usepackage[english]{babel}
\usepackage{gensymb}

\usepackage{amsfonts,amsbsy,amssymb,amsmath,graphicx,float} 
\usepackage{color}

\begin{document}

\title{Navigating Phase Space Transport with the Origin-Fate Map}

\author{Malcolm Hillebrand}
\email{hillebrand@pks.mpg.de}
\affiliation{Nonlinear Dynamics and Chaos Group, Department of Mathematics and Applied Mathematics, University of Cape Town,
Rondebosch 7701, South Africa}
\affiliation{Max Planck Institute for the Physics of Complex Systems, N\"othnitzer Stra\ss e 38, 01187 Dresden, Germany}
\affiliation{Center for Systems Biology Dresden, Pfotenhauer Stra\ss e 108, 01307 Dresden, Germany}

\author{Matthaios Katsanikas}
\affiliation{Research Center for Astronomy and Applied Mathematics, Academy of Athens, Soranou Efesiou 4, Athens, GR-11527, Greece.}
\affiliation{School of Mathematics, University of Bristol, \\ Fry Building, Woodland Road, Bristol, BS8 1UG, United Kingdom}
\author{Stephen Wiggins}
\affiliation{School of Mathematics, University of Bristol, \\ Fry Building, Woodland Road, Bristol, BS8 1UG, United Kingdom}
\affiliation{Department of Mathematics, United States Naval Academy, Chauvenet Hall, 572C Holloway Road Annapolis, MD 21402-5002, United States of America}

\author{Charalampos Skokos}
\affiliation{Nonlinear Dynamics and Chaos Group, Department of Mathematics and Applied Mathematics, University of Cape Town,
Rondebosch 7701, South Africa}

\begin{abstract}
We introduce and demonstrate the usage of the origin-fate map (OFM) as a tool for the detailed investigation of phase space transport in reactant-product type systems.
For these systems, which exhibit clearly defined start and end states, it is possible to build a comprehensive picture of the lobe dynamics by considering backward and forward integration of sets of initial conditions to index their origin and fate.
We illustrate the method and its utility in the study of a two degrees of freedom caldera potential with four exits, demonstrating that the OFM not only recapitulates results from classical manifold theory, but even provides more detailed information about complex lobe structures.
The OFM allows the detection of dynamically significant transitions caused by the creation of new lobes, and is also able to guide the prediction of the position of unstable periodic orbits (UPOs).
Further, we compute the OFM on the periodic orbit dividing surface (PODS) associated with the transition state of a caldera entrance, which allows for a powerful analysis of reactive trajectories.
The intersection of the manifolds corresponding to this UPO with other manifolds in the phase space results in the appearance of lobes on the PODS, which are directly classified by the OFM.
This allows computations of branching ratios and the exploration of a fractal cascade of lobes as the caldera is stretched, which results in fluctuations in the branching ratio and chaotic selectivity.
The OFM is found to be a simple and very useful tool with a vast range of descriptive and quantitative applications.
\end{abstract}

\maketitle


\section{Introduction}
\label{intro}
The study of phase space transport in dynamical systems has a rich and varied history.
This has been demonstrated in applications ranging from scattering dynamics in astronomy~\cite{Seoane2006,athanassoula2009rings} and beyond~\cite{Seoane2013}, to turbulent plasma flows in the context of shocks~\cite{Trotta2021} and planetary magnetosphere chorus emissions~\cite{Zonca2021}, and chemical potentials~\cite{waalkens2004direct,ezra2018sampling}.
Identifying basins of attraction in multistable systems, i.e.~the ultimate forward time fate of sets of initial conditions~\cite{Nusse1996}, has drawn attention from multiple angles~\cite{Les2022}, including considering measures of chaos based on fractality of basin boundaries \cite{Grebogi1983,Aguirre2009}, notions of basin stability from volumes of basins of attraction \cite{Menck2013}, and basin entropy as a measure of the unpredictability of the final state of a system~\cite{Daza2016,Daza2022}.
Beyond only forward time attractors, in many systems there is a clear definition of an origin state as well as the final region, where it can be useful to study this phenomenon of interstate transport, particularly in chemical reaction dynamics where reactant and product states are of note~\cite{carpenter1985,collins2014}.
In such cases, where the phase space can be partitioned into distinct regions constituting possible origin states, possible final fate states, and an intermediate transient, then the language of symbolic dynamics can be effective for describing the transitions between these states~\cite{Morse1938,LindMarcus1995}.
This allows a continuous dynamical system to be described in terms of the qualitatively significant areas of phase space, giving a practical discretisation of phase space transport by attributing separate labels to each of these regions, and allocating this label to all trajectories within the region.

In a wide variety of systems we can further take advantage of the theory of lobes which underlie long-time dynamics, by considering intersections of stable and unstable invariant manifolds corresponding to periodic orbits~\cite{Wiggins1992}.
With these manifolds partitioning the phase space into classes of trajectories having specific behaviours, this theory naturally complements the symbolic dynamics approach by proposing potentially significant regions or states to consider symbolically.
In order to compute these manifolds, one can implement classical methods using the eigenvectors of periodic orbits, or where applicable simpler techniques such as the ridges of Lagrangian descriptor gradients~\cite{mancho2013lagrangian,lopesino2017,ldbook2020}.
This manifold approach to phase space transport has for example been used to study eddies in fluid flows~\cite{Haller1998}, to compute branching ratios in chemical reactions~\cite{Katsanikas2020pre}, and explored in the context of stochastic flows~\cite{Balasuriya2018}.

While applicable to a diverse array of physical problems, the modelling of chemical reactions is a prime example of the utility of these transport analyses, particularly due to transition state theory~\cite{wigner1938,waalkens2007}.
The mathematical description of this transition state leads to a dividing surface, or surface of no return, which can be directly constructed in some models~\cite{waalkens2004direct}.
Applying these techniques to a variety of potential energy surfaces has led to an increased understanding of chemical dynamics~\cite{wiggins2016,Agaoglou2019,agaoglou2020phase,katsanikas2020a,Katsanikas2023}.

In this work, we introduce the notion of the origin-fate map (OFM) as a direct method for exploring phase space transport, validate it by comparing to previous results, and explore some of its capabilities in the context of a simple chemical reaction potential model.
Section~\ref{sec:ofm} defines the OFM and how it is computed, and in Section~\ref{sec:transport_in_a_symmetric_caldera} we show that the method can be used to reproduce and extend previous findings in a caldera potential.
Section~\ref{sec:transport_from_the_periodic_orbit_dividing_surface} describes the advantage of considering the map on a dividing surface of the reaction potential, and provides further illustration of the qualitative and quantitative insights it provides.
Finally, Section~\ref{sec:summary_and_conclusions} concludes the work by summarising the method and some of its use cases.

\section{The Origin-Fate Map}
\label{sec:ofm}
An inherent difficulty with any manifold-based transport characterisation technique is the problem of correctly classifying the various manifold lobes.
Similarly, in situations where the manifolds may be complex and difficult to extract, classifying regions of phase space based on their transport characteristics is challenging.
What we propose here is a simple method to identify asymptotic behaviours of trajectories through an origin-fate map for reactant-product type systems, extending ideas from both lobe dynamics and studies of basins of attraction.
Fundamentally, for a given initial condition, typically on a surface of section, we can identify an \textit{origin}, associated with an initial reactant region, for instance an attractor, a potential well, or entry channel, and a \textit{fate}, associated with the final product region or exit channel.
We can directly identify the origin through backward integration, and the fate via forward integration.
We then associate an \textit{origin-fate index} with each of these points on the surface of section, consisting of two numbers indexing the origin and fate channels of the trajectory.
We note that using forward integration to investigate the fates of trajectories is well explored, particularly through basins of attraction (see for instance~\cite{Nusse1996,Seoane2006,carpenter2018}), but to the best of our knowledge no comprehensive, detailed methods have been proposed to deal with both origins and fates of trajectories.

\begin{figure}[tb]
    \centering
    \includegraphics[width=\columnwidth]{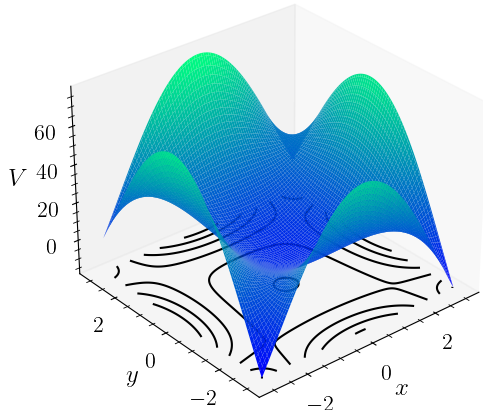}
    \caption{The unstretched ($\lambda=1$) caldera potential energy surface $V(x,y)$ from Eq.~\eqref{eq:caldera}, showing the four exit channels at each corner and the central well. The contours of the potential are shown in black below the 3D surface.}
    \label{fig:3dpotential}
\end{figure}

For simplicity and concreteness, we present and describe the method in the context of a specific Hamiltonian caldera model satisfying the reactant-product criterion.
This caldera potential energy surface takes motivation from organic chemical reactions, where it has been seen to describe several rearrangements of molecules~\cite{gold1988,doubleday1997,reyes2002,doubleday2006}.
The notions discussed here are straightforwardly extensible to different open Hamiltonian models with escapes, and to dissipative systems with forward- and backward-time attractors; we emphasise that the technique works for Hamiltonian and non-Hamiltonian systems alike.
The Hamiltonian of this caldera system is given in terms of a unit mass test particle with position $(x,y)$ and momentum $(p_x,p_y)$ by~\cite{collins2014}
\begin{align}
    \label{eq:caldera}
    \mathcal{H}(x,y,p_x,p_y) &= \frac{p_x^2}{2} + \frac{p_y ^2}{2} + V(x,y); \\
    V(x,y) &= c_1(\lambda^2 x^2+y^2) +c_2y \nonumber\\
    &-c_3(\lambda^4x^4+y^4-6\lambda^2x^2y^2). \nonumber
\end{align}
In Eq.~\eqref{eq:caldera}, $V(x,y)$ describes a potential energy surface with two higher-energy index 1 saddles corresponding to entrances to the caldera, a central minimum, and two lower energy index 1 saddles corresponding to exits.
The value of the Hamiltonian is the overall energy of the system $E$.
The parameter $\lambda$ controls the stretching of the surface, with $\lambda=1$ giving the unstretched potential, and decreasing values of $\lambda$ resulting in increased stretching of the potential energy surface.
To maintain the four-escape structure, the parameters are taken as $c_1 = 5,\ c_2 = 3,\ c_3 = -0.3$ as in Refs.~\cite{katsanikas2018,katsanikas2019}.
Chemically relevant trajectories enter and exit through one of the four corners defined by the saddles and bounded by isoenergetic surfaces.
The only alternative is trajectories which are trapped within the central basin of the caldera in either forward or backward time.
The unstretched potential energy surface is shown in Fig.~\ref{fig:3dpotential}, demonstrating these features: The four escape channels, the central well, and the $y-$asymmetry creating the higher and lower energy channels. 

\begin{figure}[tb]
    \centering
    \includegraphics[width=.8\columnwidth]{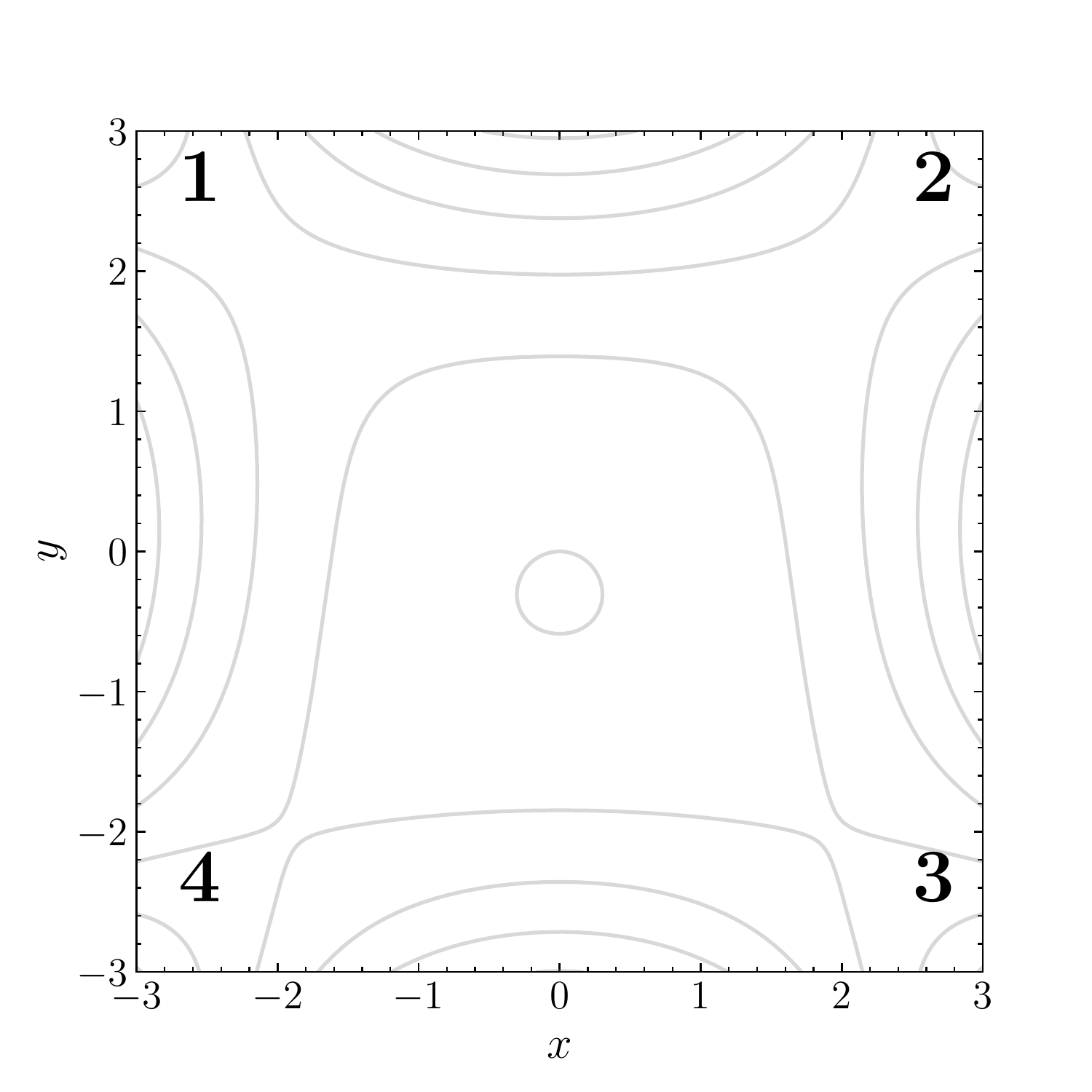}
    \caption{The contours of the unstretched ($\lambda=1$) caldera potential~\eqref{eq:caldera} as in Fig.~\ref{fig:3dpotential}, marked with the origin-fate indices. So the top left channel is indexed as channel 1, the top right indexed 2, the bottom right indexed 3, and the bottom left is indexed 4.}
    \label{fig:entrylabels}
\end{figure}

In the caldera, we assign origin/fate indices to the four channels as per Fig.~\ref{fig:entrylabels}, moving clockwise from top left.
The top left is indexed as 1, where we have the top left saddle (TLS), the top right indexed 2 (TRS), the bottom right indexed 3 (BRS), and the bottom left indexed 4 (BLS).
For example, if we had a trajectory that enters through the top right channel, and then exits through the bottom left, it would be indexed as $2-4$: It has origin channel 2, and fate channel 4.
We note that entering or leaving through one of these channels corresponds to a particular symbolic dynamics state, and consequently the origin-fate index essentially describes a symbolic trajectory.

This allows us to index each possible combination of origin-fate channels using 16 identifiers, with one additional identifier for the case where (up to the finite integration time) the trajectory remains trapped in the caldera in either the forward or backward integration.
By assigning colours to each of these identifiers, it is possible to plot the surface of section coloured by the origin-fate index of each point.
While this explanation refers to the specific case of the caldera with four entry and exit channels, clearly similar descriptions would hold for a general system with $n$ possible origins and $k$ possible fates, requiring $n\times k$ colours.

The OFM algorithm for a given model with set parameters can be summarised in the following steps.
\begin{enumerate}
    \item Choose a surface on which to take a discrete grid of initial conditions.
    \item For each initial condition
    \begin{enumerate}
        \item Integrate forward until it reaches an ending fate state, or reaches the cutoff time.
        \item Integrate backward until it reaches a starting origin state, or reaches the cutoff time.
        \item Record the origin-fate index based on these results.
    \end{enumerate}
    \item For low-dimensional visualisation: Colour the surface based on the origin-fate indices.
    \item For quantification of trajectory behaviours: Compute the fraction of initial conditions with each origin-fate index.
\end{enumerate}

The definition of an origin/fate state will need to be decided based on the properties of the system, and whether this state is an attractor, escape channel, spatially localised region etc.

\section{Transport in a Symmetric Caldera} 
\label{sec:transport_in_a_symmetric_caldera}
As a combined validation, demonstration and explanation of the OFM, we consider the caldera potential~\eqref{eq:caldera} which has been explored extensively as a transition-state model for chemical reactions with two reactant states and two product states in a symmetric~\cite{katsanikas2018,katsanikas2020b} and asymmetric~\cite{Katsanikas2022} formulation.
In the unstretched caldera with $\lambda=1$, the phenomenon of ``dynamical matching'' is observed, which means that all reactive trajectories cross the caldera diagonally -- in the notation of origin-fate indices, the only possible reactive cases are $1-3$ and $2-4$ (cf.~the labels in Fig~\ref{fig:entrylabels}).
It has been shown that as this potential surface is stretched by decreasing $\lambda$, lobes are formed by the unstable manifold of the unstable periodic orbit (UPO) corresponding to the TRS and stable manifolds from periodic orbits in the central minimum~\cite{katsanikas2018,katsanikas2020b}.
The formation of these lobes result in the breaking of dynamical matching, the appearance of trajectories that enter from the top right (left) and leave through the bottom right (left) -- i.e.~allowing $1-4$ and $2-3$ origin-fate indices.
Here, we demonstrate that the OFM shows this result directly, while additionally providing a more detailed description of the lobe structures responsible for phase space transport.

\begin{figure}[tb]
    \centering
    \includegraphics[width=\columnwidth]{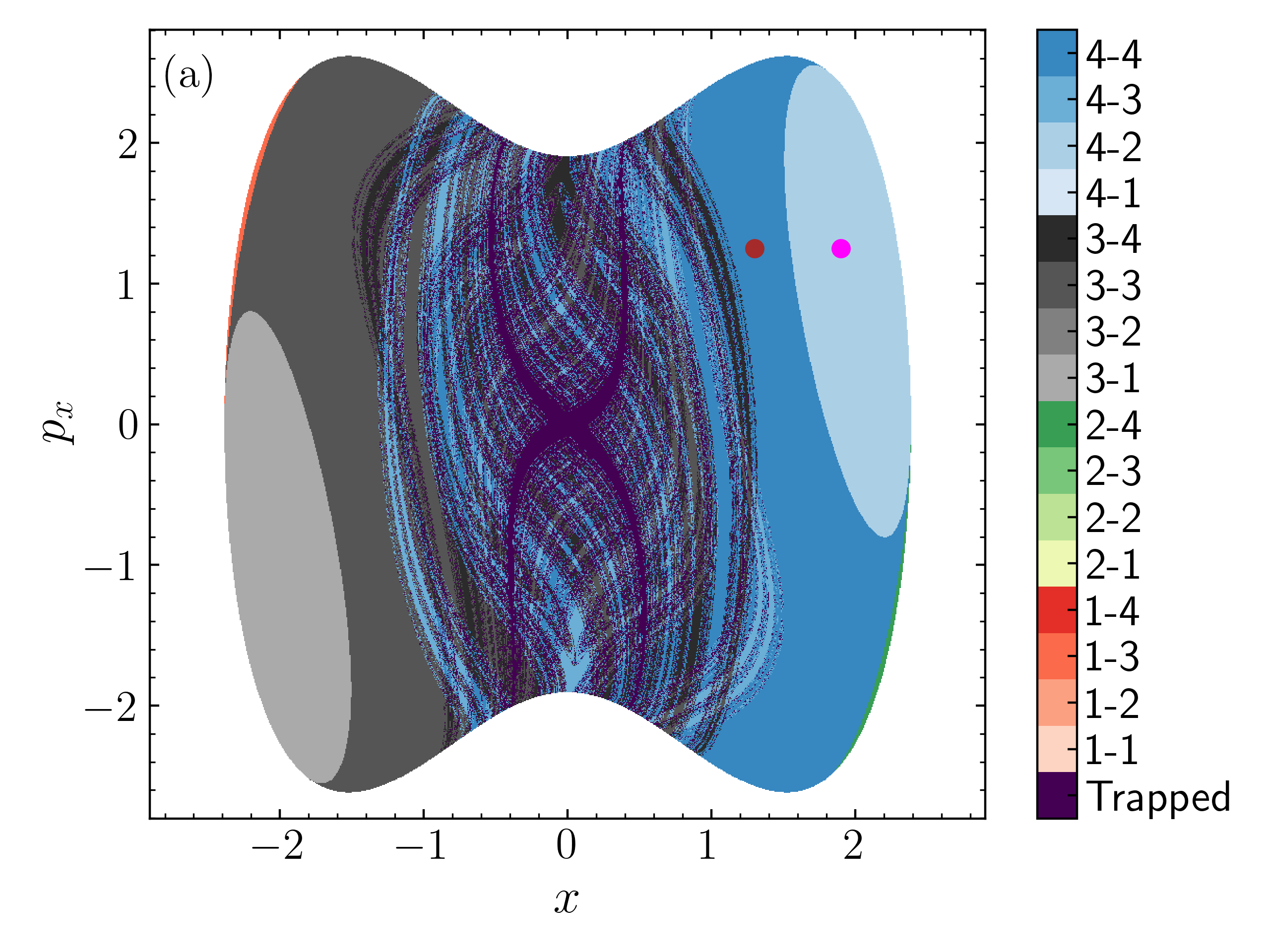}
    \includegraphics[width=0.8\columnwidth]{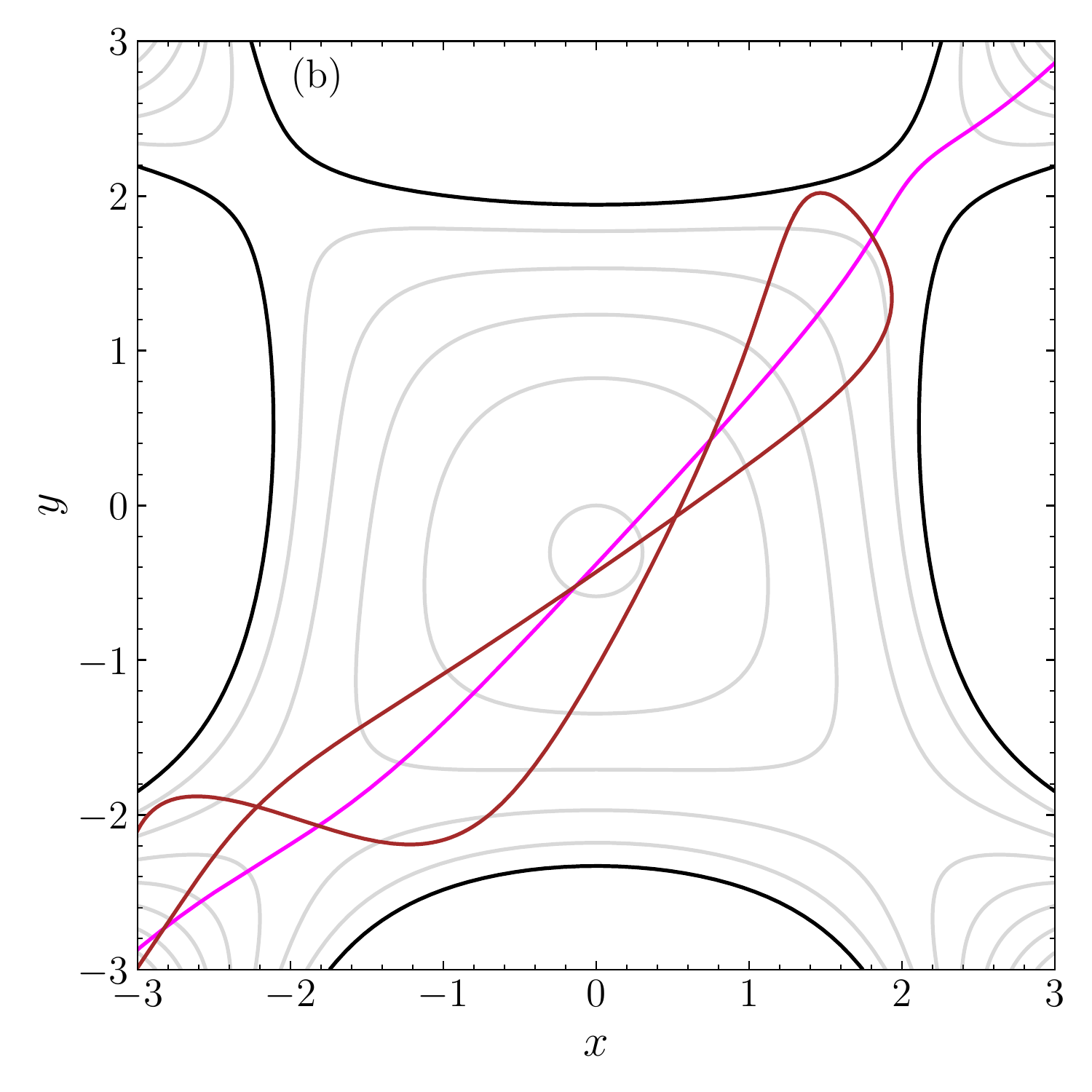}
    \caption{(a) The OFM using an integration time of $\tau=20$ for the Poincar\'e surface of section at $y=1.88409,\ p_y>0$ of the unstretched ($\lambda=1$) caldera for energy $E=29$ showing a complex array of structures, with the positive $p_y$ resulting in a majority of trajectories entering from below (indices 3 and 4). (b) Two representative trajectories in configuration space, showing the two points marked on the surface of section in (a). Black lines correspond to the energy boundary.}
    \label{fig:SC_full}
\end{figure}
\begin{figure}[tb]
    \centering
    \includegraphics[width=\columnwidth]{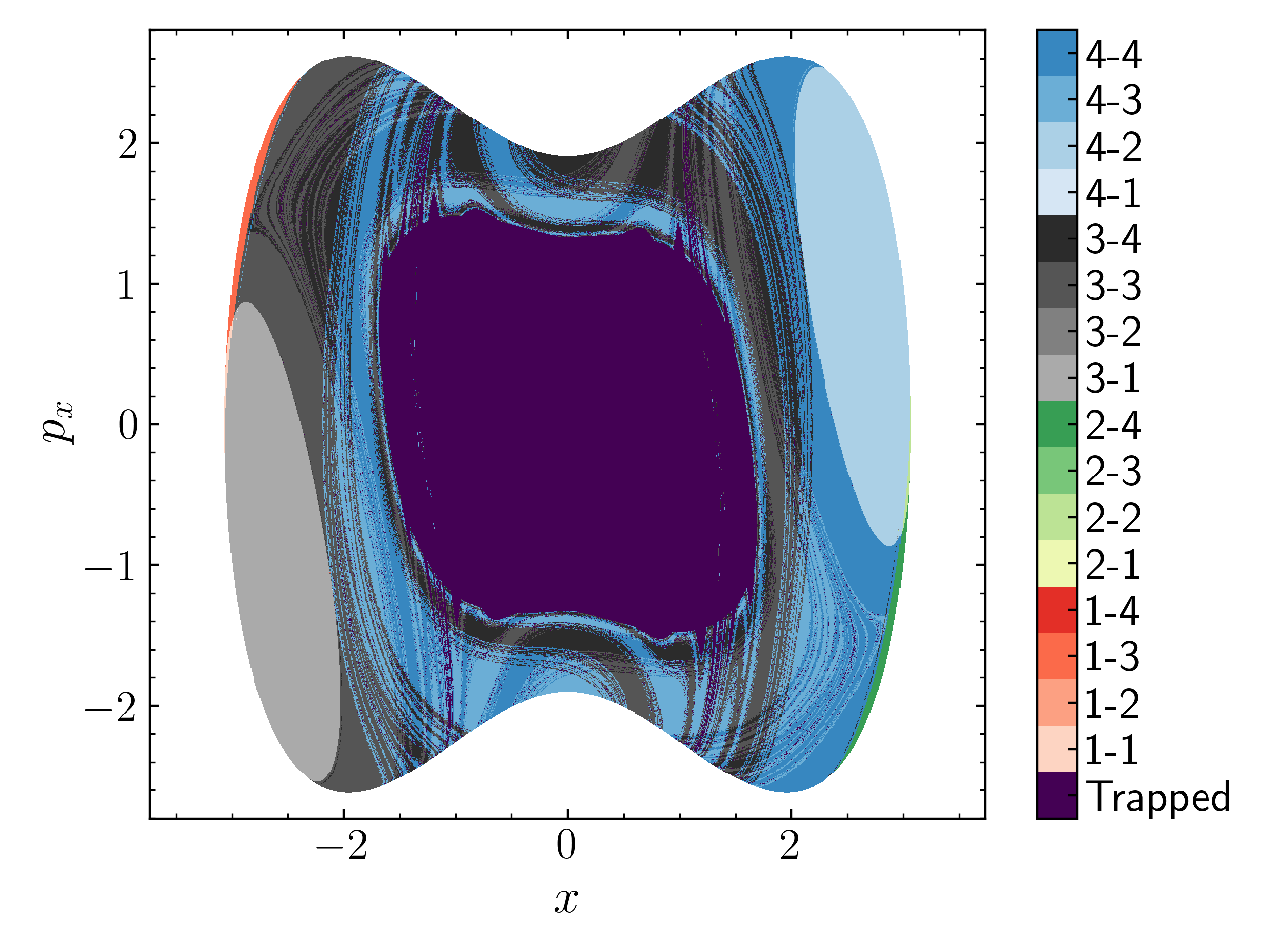}
    \caption{Similar to Fig.~\ref{fig:SC_full}, but for the critical stretching of the potential, $\lambda=0.778$, where the first signs of breaking dynamical matching appear.}
    \label{fig:SC_fullstretch}
\end{figure}

In Fig.~\ref{fig:SC_full}(a), we see a visualisation of the OFM for the caldera potential with energy $E=29$, on the surface of section defined by $y=1.88409,\ p_y>0$ (after~\cite{katsanikas2020b}).
A maximum integration time of $\tau=20$ was used, i.e.~integrating until $t=\tau$ or the trajectory exits the caldera region.
Firstly, this shows how the OFM clearly delineates distinct behaviours, with the colour of an initial condition corresponding to its origin [see the colour bar in Fig.~\ref{fig:SC_full}(a)], and the gradient of the intensity denoting the fate.
Distinct lobes are visible, with simple structures on the edges of the filled space and more complex entanglements in the central region.
Figure~\ref{fig:SC_full}(b) shows two of the distinct origin-fate behaviours by following the evolution of the initial conditions through configuration space.
These initial conditions are indicated by the coloured points on the surface of section in Fig.~\ref{fig:SC_full}(a).
Essentially, this encapsulates the utility of the OFM -- instead of requiring exhaustive trial and error tests of individual initial conditions, we can clearly visualise the pertinent behaviours in a symbolic dynamics sense.

Further, the transport characteristics of this caldera potential have already been studied through manifold dynamics, with both classical methods~\cite{katsanikas2018}, and Lagrangian descriptors~\cite{katsanikas2020b}.
In order to provide a clear comparison with these previous results, and to demonstrate the practicality of the OFM in adding to the existing tools, we consider the stretching of the potential surface by decreasing the value of $\lambda$ in Eq.~\eqref{eq:caldera}.
When $\lambda$ attains its critical value of $\lambda\approx 0.778$, it has been observed that the unstable manifold of the UPO corresponding to the TRS is intersected by manifolds from periodic orbits in the central region, resulting in a breaking of dynamical matching (i.e.~reactive trajectories no longer conform solely to $1-3$ or $2-4$ classes)~\cite{katsanikas2018,katsanikas2020b}.
In Fig.~\ref{fig:SC_fullstretch} we show the OFM for the same surface of section as in Fig.~\ref{fig:SC_full}, at this critical value.
There are several features that are immediately apparent -- the greater proportion of trapped trajectories (purple regions) for this finite time ($\tau=20$), and the expansion of the ``complex central region'', as manifolds from periodic orbits in the center reach further in the $x-$direction.
We note that a longer time can be used to reduce the number of trapped trajectories, or the forward/backward nature of the trapping could also be directly tracked if this is of interest.

\begin{figure*}[tb]
    \centering
    \includegraphics[width=0.9\columnwidth]{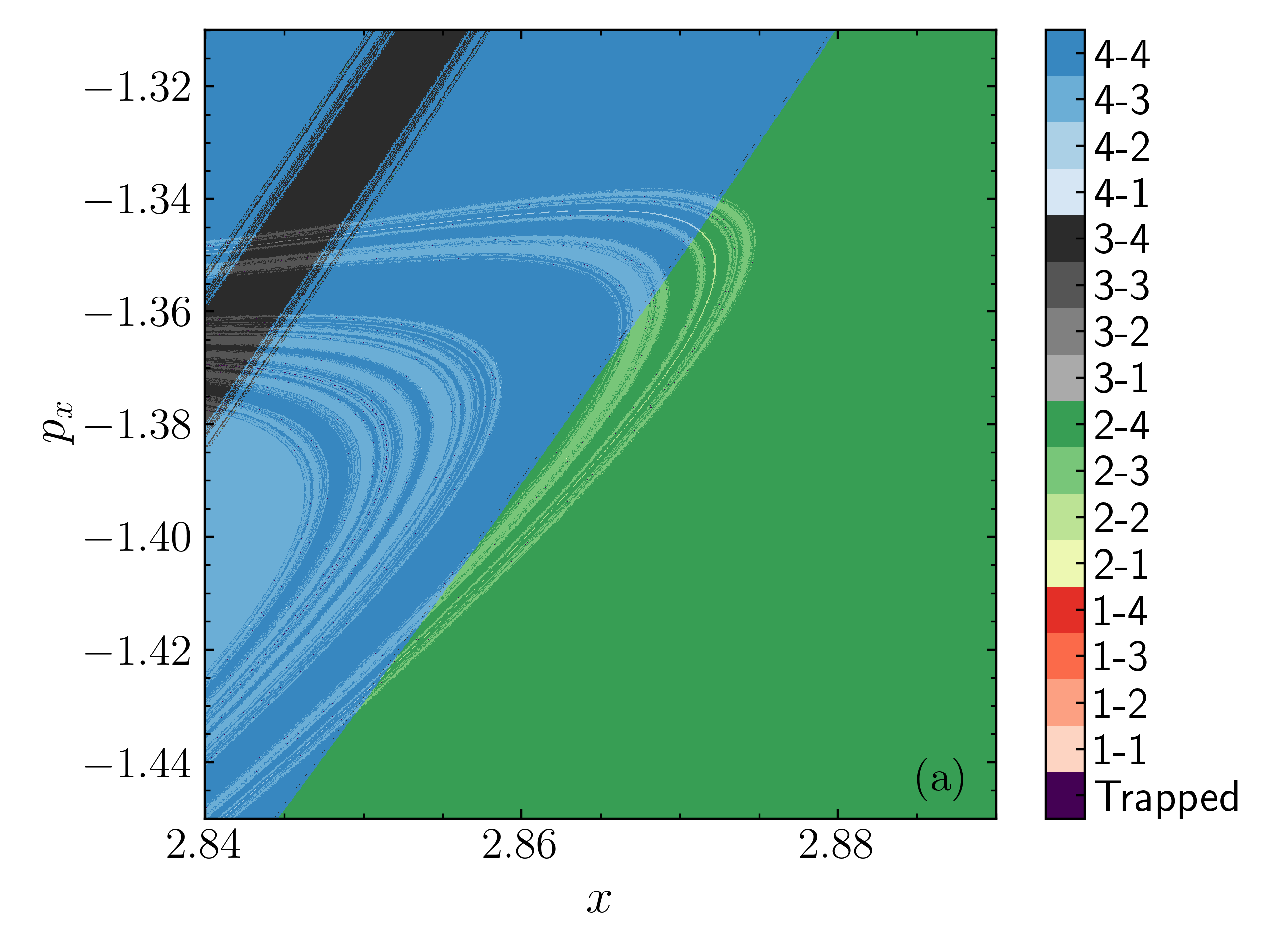}
    \includegraphics[width=0.9\columnwidth]{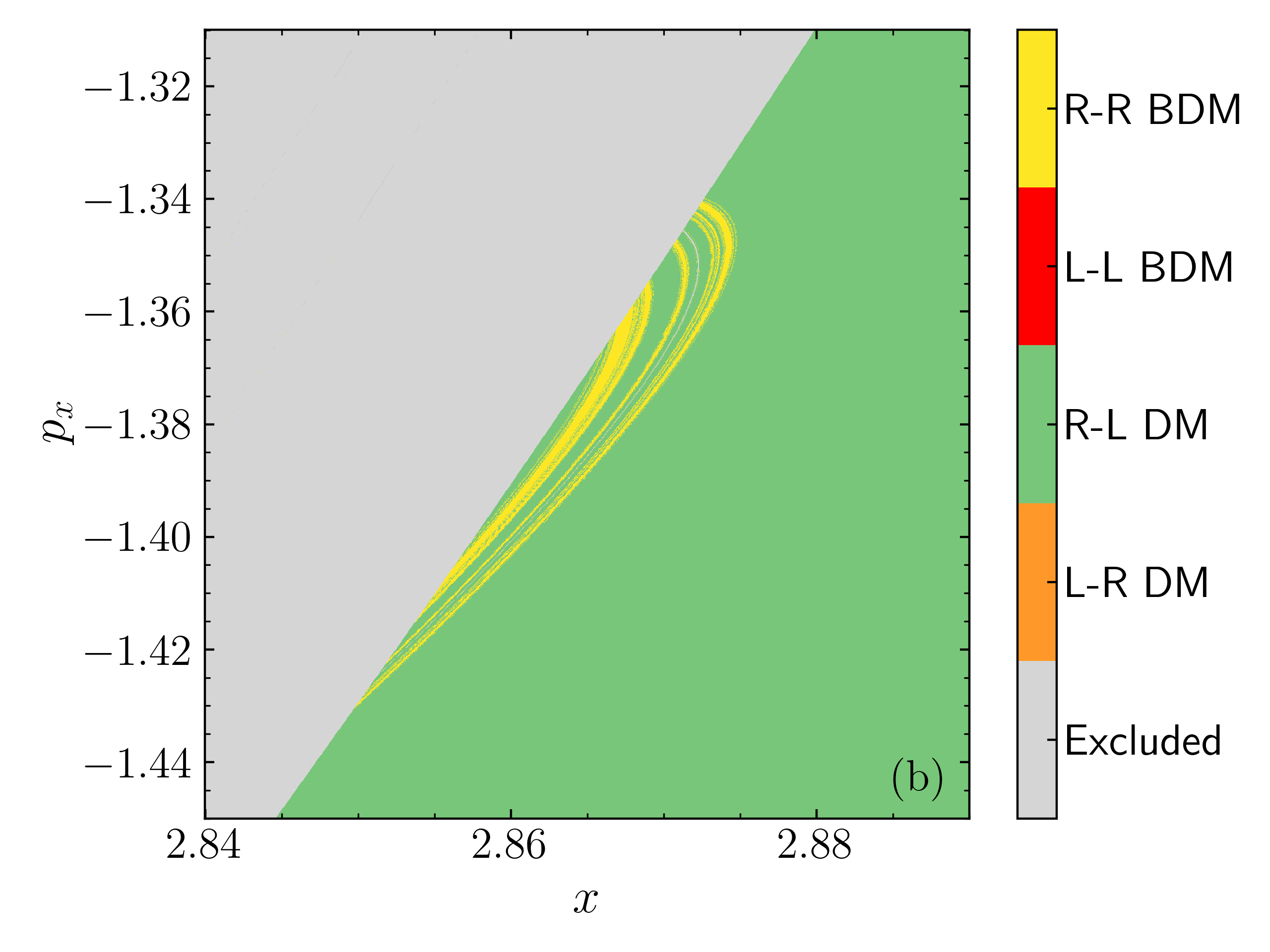}
    \includegraphics[width=0.9\columnwidth]{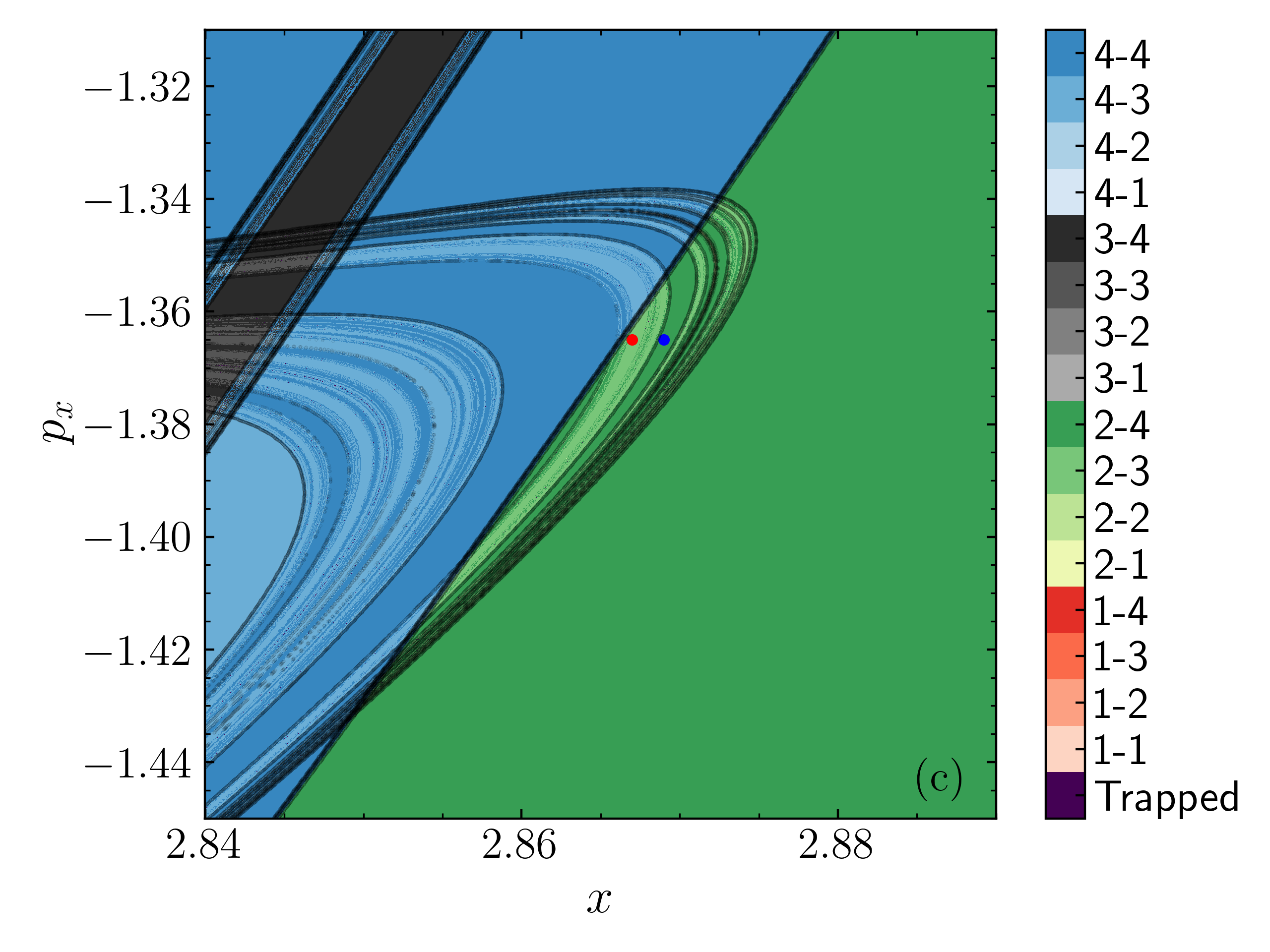}
    \includegraphics[width=0.9\columnwidth]{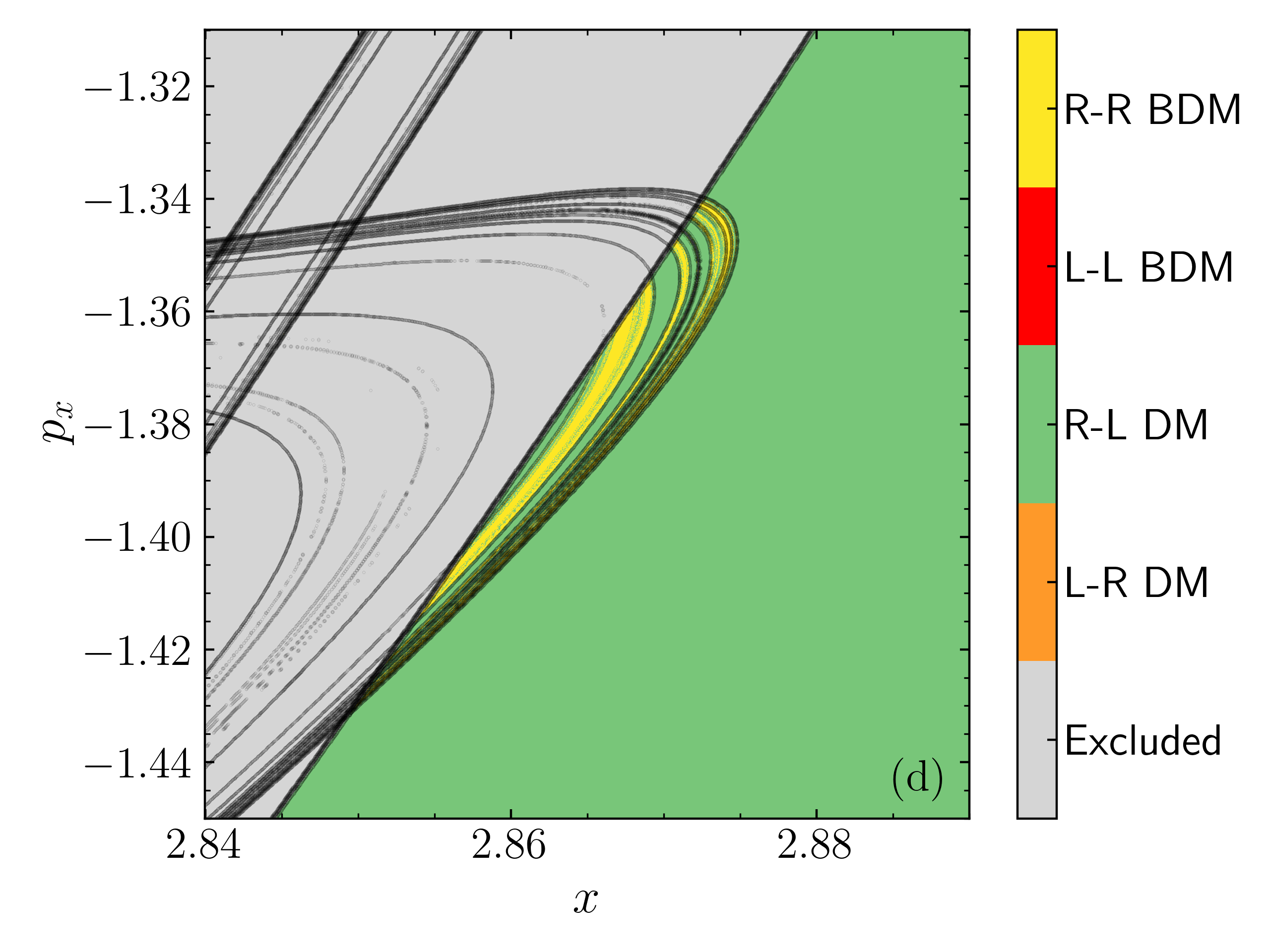}

    \caption{The OFM for the region of the Poincar\'e surface of section in Fig.~\ref{fig:SC_fullstretch} where the lobes appear. (a) The full map, showing all origin-fate indices. (b) The COFM, where we reduce the colors to behaviors of interest -- Right-Right broken dynamical matching (yellow), right-left dynamical matching (green), left-left broken dynamical matching (red) and left-right dynamical matching (orange). Other behaviours are colored grey. (c) The same as panel (a), but with the manifolds overlaid as black lines, and two specific trajectories marked as red and blue points (see Fig.~\ref{fig:trajs}). (d) The same as panel (b) with the manifolds overlaid, demonstrating the precise correspondence between the manifold lobes and the different transport behaviours.}
    \label{fig:SC_lobes}
\end{figure*}

\begin{figure}[tb]
    \centering
    \includegraphics[width=0.9\columnwidth]{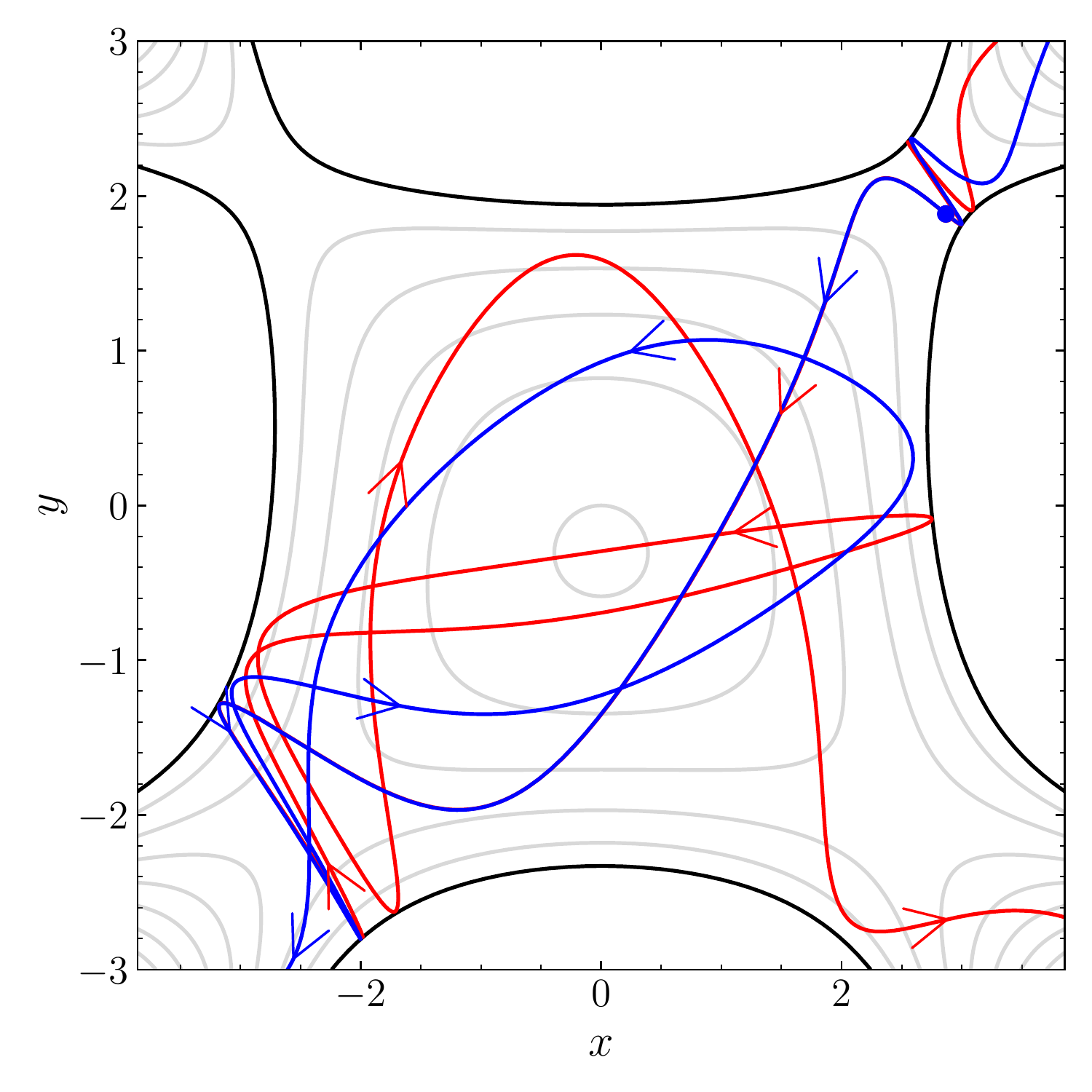}
    \caption{The configuration space projection of two trajectories that cross the surface of section shown in Fig.~\ref{fig:SC_lobes} in the lobes formed by the unstable invariant manifold corresponding to the UPO of the TRS and the manifolds from the central area. The contours of the potential are shown in grey, with the black lines demarcating the energy boundary. The red trajectory corresponds to a point in a $2-3$ lobe, and the blue trajectory to a point in the $2-4$ lobe. The initial conditions are marked by (overlapping) red and blue dots. Note that the initial part of the forward evolution also overlaps completely on this plot. The arrows indicate the forward-time direction of the motion.}
    \label{fig:trajs}
\end{figure}

What is only slightly apparent at this scale is the expected intersection of manifolds resulting in the breaking of dynamical matching, and the creation of $2-3$ (lighter green colour) or $1-4$ (deep red colour) origin-fate classes.
In order to properly investigate this, we consider a zoomed in region of the bottom right corner of the surface of section in Fig.~\ref{fig:SC_fullstretch}, near the small green $2-4$ lobe.
Figure~\ref{fig:SC_lobes}(a) shows the OFM of this zoomed in region, showing exactly the manifold-predicted behaviour: There are lobes of broken dynamical matching forming, with a $2-3$ origin-fate index.
Here we can directly observe that these lobes are complex, as there are rapidly alternating origin-fate regions described by manifolds of higher-order periodic orbits.
In fact, as will be discussed later, there is a chaotic fractal structure of manifolds in this region, resulting in unpredictable transport effects.

Since we are primarily concerned only with dynamical matching and the breaking thereof, a natural way to adapt the OFM is to colour initial conditions on the surface based on these classifications.
Concretely, we will consider only four relevant, or ``active'' cases: Right-left dynamical matching (R-L DM, $2-4$), left-right dynamical matching (L-R DM, $1-3$), right-right broken dynamical matching (R-R BDM, $2-3$), and left-left broken dynamical matching (L-L BDM, $1-4$).
All other behaviours will be ignored by simply colouring those initial conditions grey.
Fig.~\ref{fig:SC_lobes}(b) shows this classification OFM (COFM), where the behaviour of most interest, the R-R BDM, is shown in yellow.
Here our attention is clearly drawn to the key characteristics of the complex lobes.
While exactly the same information is contained here as in the full OFM in Fig.~\ref{fig:SC_lobes}(a), the removal of extraneous distracting colours makes the results more clearly comprehensible, and this is likely the best path for visualisations when there are \textit{a priori} behaviours of particular interest.

Crucially, we are able to show that the OFM correctly reproduces, and in fact enhances, the previous findings using manifolds.
In Figs.~\ref{fig:SC_lobes}(c) and (d) we see the OFM and COFM now overlaid with the manifolds computed through Lagrangian descriptor gradients, demonstrating that the coloured lobes of the OFM are delineated by these manifolds.
The OFM directly provides the classification of each of these lobes, without the need for trial and error.
Further, it reveals the complex sub-lobe structures that require computation of higher-order manifolds to see through classical methods, allowing us to understand the detailed dynamics without difficult manifold identification.

To illustrate the complex and sensitive dynamics, two particular initial conditions in different lobes are taken from the OFM of the caldera at the critical stretching point $\lambda=0.778$, marked by red and blue points respectively in Fig.~\ref{fig:SC_lobes}(c), and their full configuration space evolutions are shown in Fig.~\ref{fig:trajs} by going forward and backward in time.
The two initial conditions are taken as $(x,y,p_x,p_y) = (2.867,1.88409,-1.365,0.86563646)$ for the $2-3$ case (red), and $(2.869,1.88409,-1.365,0.8527248)$ for the $2-4$ case (blue).
While the trajectories are initially almost identical after crossing the surface of section at $y=1.88409$, they eventually diverge completely and exit through different channels, as indicated by their respective lobes in the OFM.
This again emphasises a principal use case for the OFM -- essentially automating the identification of these very narrow sub-lobes of different behaviours.

Another important note is the identification of fractal regions in the transport lobes.
If we consider only the trajectories emanating from the top right (i.e.~the green right hand side in Fig.~\ref{fig:SC_lobes}), then the OFM simply describes the exit channels as attractors, and we can directly identify the regions of infinitely dense manifolds as fractal boundaries to the basins of attraction of these channels.
These fractal boundaries are evident from the interleaved yellow and green regions of Fig.~\ref{fig:SC_lobes}(b), where the yellow colour corresponds effectively to the basin of attraction of the bottom right, and the green to the basin of the bottom left.
Consequently, we are able to observe that even in this simple model, there is a chaotic element present when attempting to predict the ultimate fate of particular initial conditions. 
The numerical confirmation of this sensitive nature is much more directly feasible through the OFM than any manifold-based attempts, where it becomes progressively more difficult to find or compute the higher-order manifolds responsible for the dense fractal structure.

\section{Transport from the Periodic Orbit Dividing Surface} 
\label{sec:transport_from_the_periodic_orbit_dividing_surface}

A significant decision in any method of classifying transport dynamics or basins of attraction is how to choose the (hyper)surface on which to classify initial conditions.
As in the above discussion in Section~\ref{sec:transport_in_a_symmetric_caldera}, often a well-chosen phase space slice at a particular coordinate yields an effective Poincar\'e surface of section.
However, the specific context of chemical reaction theory suggests a natural choice for the surface -- the dividing surface, which exists at the transition state, separating ``reactant'' phase space from ``product'' phase space.
As this dividing surface is a surface of no return, which satisfies a local no-recrossing property, once a trajectory crosses the surface it cannot exit via the channel it entered from without an extended journey through the phase space.

In the particular context of potentials with index 1 saddles (such as the caldera considered here), we are able to construct the dividing surface in phase space through finding UPOs corresponding to each saddle, known as periodic orbit dividing surfaces (PODS)~\cite{wiggins2001}.
In configuration space, the projection of these UPOs provides a direct dividing surface between reactants and products~\cite{Pechukas73,Pollak78,Pechukas79}.
In phase space, the dividing surface of a particular UPO is given by the constant-energy hypersurface in phase space bounded by the position coordinates of the UPO.
Explicitly, to construct the PODS we can consider every point along the UPO in configuration space, and sweep out every energetically permitted set of momentum values at this point.
The union of all these phase space points constitutes the PODS.

In principle, this PODS can then be used as a surface of section in the phase space, with notable benefits.
In particular, if we wish to study all reactive trajectories emanating from, or exiting through, a certain channel, we can compute the OFM on the corresponding dividing surface and we are guaranteed to catch every relevant trajectory.
This provides an unambiguous context for quantifying branching ratios, or the fraction of trajectories that correspond to certain origin-fate behaviours, by giving a set of initial conditions that include only reactive trajectories, and in fact include \textit{all} reactive trajectories~\cite{deAlmeida1990,Ampawan1993,Wadi1997}

For the caldera, we can compute the PODS corresponding to each saddle, and extend the analysis performed on the constant $y$ surface of section.
In practise, this consists of two steps -- firstly identifying and finding the UPO, and then compute origin-fate indices for every initial condition lying on the dividing surface.

\begin{figure}[tb]
    \centering
    \includegraphics[width=\columnwidth]{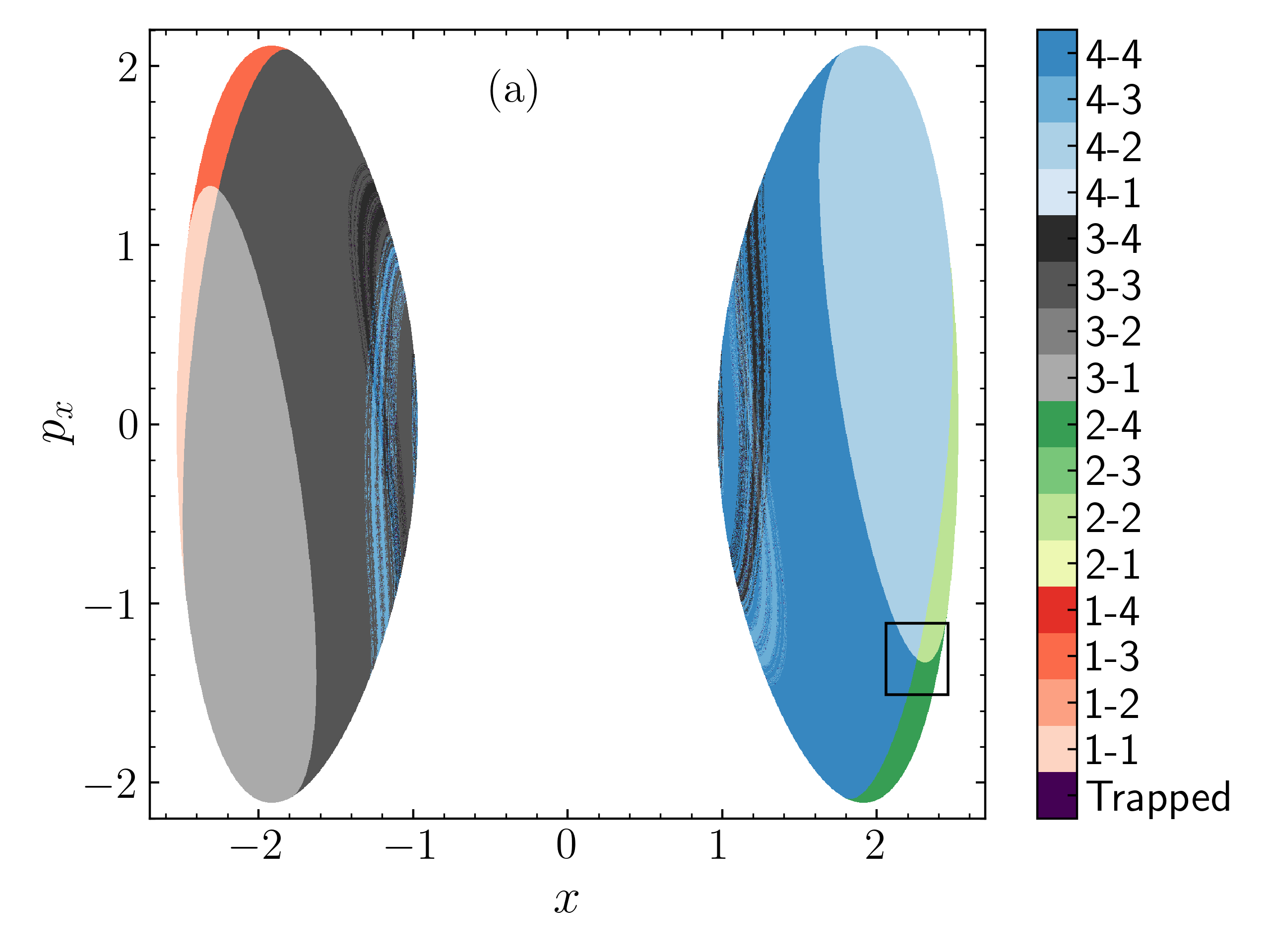}
    \includegraphics[width=\columnwidth]{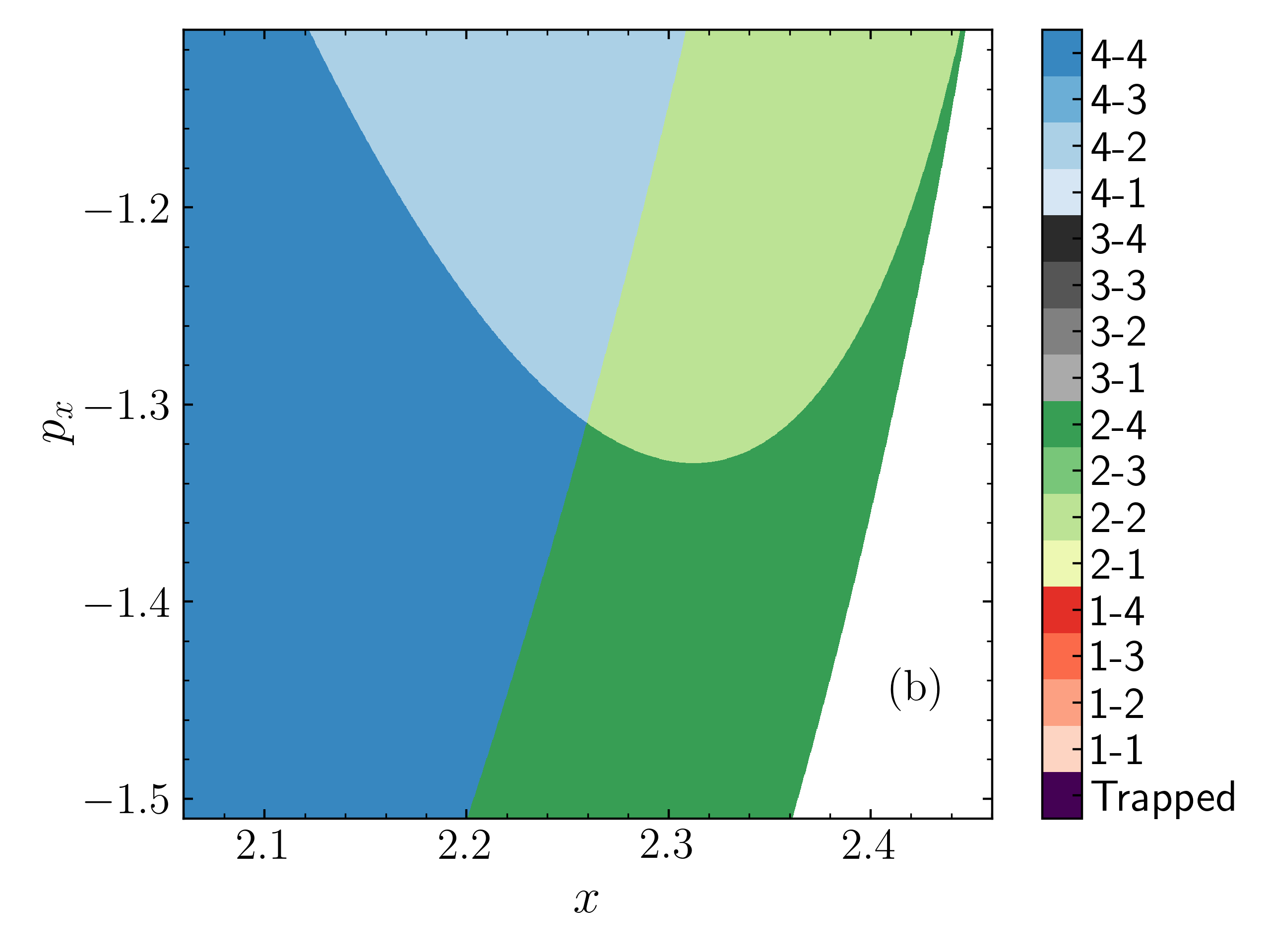}
    \caption{(a) The OFM for the unstretched caldera with $E=29$, on the surface of section $y=2, p_y>0$. (b) An enlargement of the rectangular region marked in panel (a), showing the four distinct quadrants of origin-fate behaviours characteristic of an unstable periodic orbit (see text).}
    \label{fig:UPO_identify}
\end{figure}

\begin{figure*}[tb]
    \centering
    \includegraphics[width=0.45\textwidth]{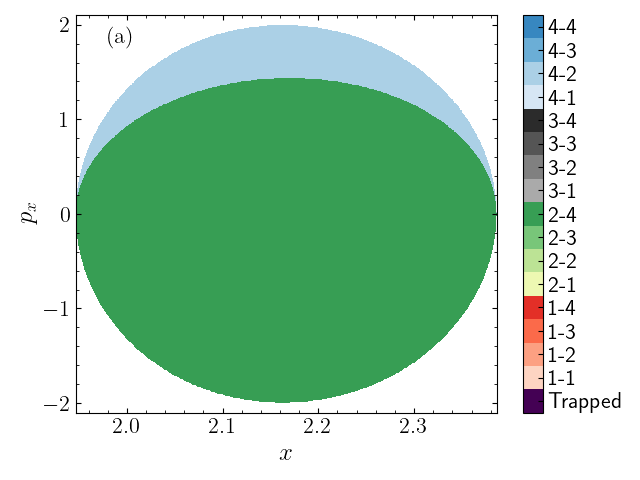}
    \includegraphics[width=0.45\textwidth]{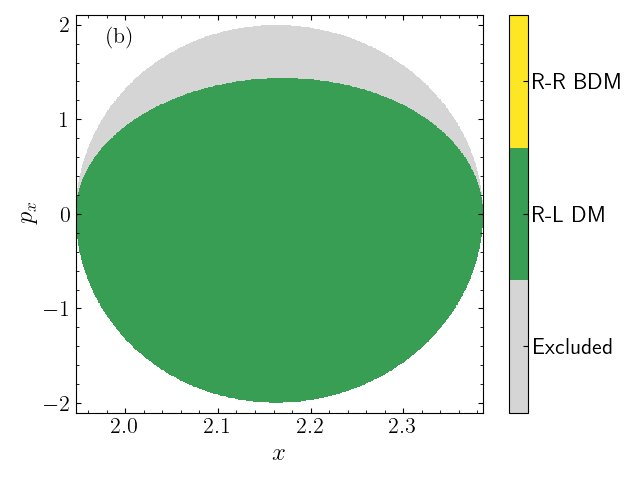}
    \includegraphics[width=0.45\textwidth]{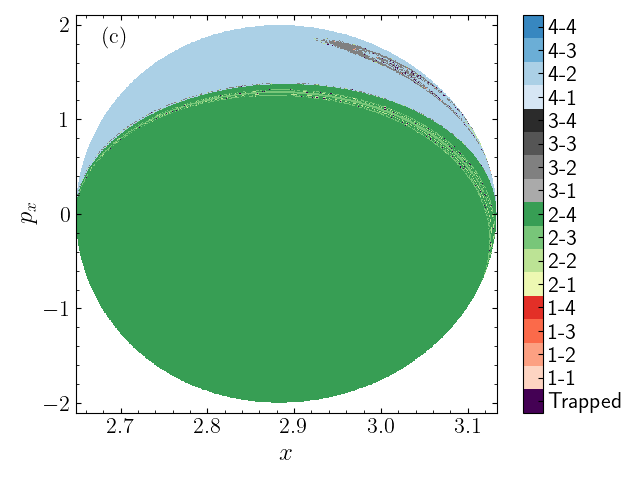}
    \includegraphics[width=0.45\textwidth]{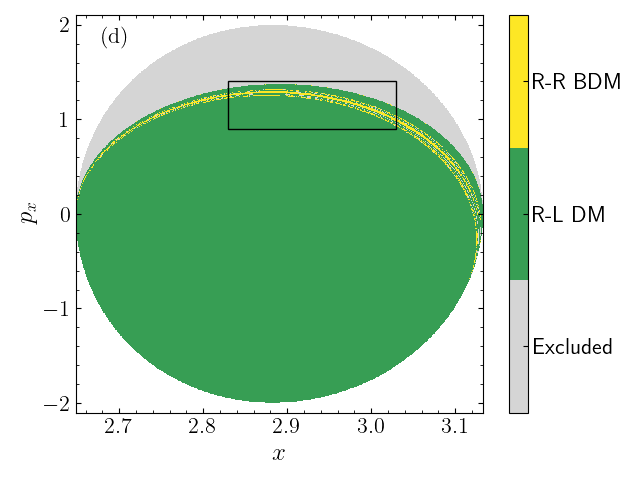}
    \caption{The $x-p_x$ projection of the PODS associated with the top right saddle with $E=29$ and $p_y<0$, coloured according to the origin-fate map. Panels (a) and (c) show the full map, with $\lambda=1$ and $\lambda=0.75$ respectively. (b) and (d) show the classification map for these $\lambda$ values, restricting options to trajectories across the caldera, 2-4 (R-L DM), and trajectories breaking dynamical matching, and exiting to the right 2-3 (R-R BDM). The rectangle in (d) demarcates the region shown in more detail in Fig.~\ref{fig:SC_zoom}}
    \label{fig:PODS}
\end{figure*}

\subsection{Finding Periodic Orbits} 
\label{sub:finding_periodic_orbits}

The first step is aided significantly by studying the OFM on a constant-coordinate surface of section.
It is known that the stable and unstable invariant manifolds of a given UPO partition the dynamics~\cite{Wiggins1992}, and the OFM gives precisely this asymptotic symbolic partition (i.e.~the origin-fate index).
Thus, by inspection of the OFM on a Poincar\'e surface of section, we can in fact identify potential UPOs which intersect the surface.
If there is a single intersection of the UPO with the surface of section, this intersection will appear as a ``corner node'', a point where the intersection of stable and unstable manifolds results in four distinct quadrants of origin-fate behaviours, corresponding to the particular UPO.
Here it is important to recall that crossing an unstable invariant manifold results in a change in the \textit{fate} index, as the unstable manifolds govern forward-time transport, and correspondingly crossing a stable invariant manifold results in a change in the \textit{origin} index.
For example, in the caldera the UPO corresponding to the TRS (recall this origin-fate channel is indexed 2) can be identified by taking a surface of section and finding such a corner node where two neighbouring quadrants have origin index 2, including one quadrant that also has fate index 2.
The four quadrants, as one rotates around the UPO intersection point, will thus be indexed $2-X$, $X-X$, $X-2$, $2-2$, where $X$ represents the index of the other channel.
This is demonstrated in the OFM segment on the surface of section $y=2, p_y>0$, Fig.~\ref{fig:UPO_identify}, where we see that due to dynamical matching the two channels involved here are 2 and 4.
The colouring technique of the OFM makes the quadrant matching in Fig.~\ref{fig:UPO_identify}(b) clear, since by rotating around the corner node, either the colour or the gradient changes, but not both.
These correspond to fate changes (colour), or origin changes (gradient), from the manifold crossing.
Thus we see if we start in the lower right quadrant, and rotate clockwise, we have the sequence described above, $2-4$, $4-4$, $4-2$, $2-2$.
This plot amply demonstrates the point that we are able to make a very accurate initial guess for the position of the UPO on this surface of section, given by the corner node of this colouring, which then enables a root-finding approach to precisely locate the UPO.

\subsection{Using the OFM on a Dividing Surface} 
\label{sub:the_ofm_on_a_pods}

Once the UPO has been found, the computation of the OFM is relatively straightforward in two degrees of freedom.
In this model, the coordinates of the UPO give a fixed relationship between $x$ and $y$ (in practice we have a list of coordinate pairs from our UPO algorithm), defining a surface of section.
With this we can follow the algorithm described above, where for each of these $(x,y)$ pairs we can compute the origin-fate indices on a grid of energetically allowed momenta.
This procedure yields the OFM on a surface, in exactly the same way as we would have at a constant coordinate.

In Fig.~\ref{fig:PODS} we show the OFM computed on the PODS corresponding to the TRS, with $E=29$, $p_y<0$, and presented in the $x-p_x$ projection (recall that $y$ is defined by $x$).
Figures~\ref{fig:PODS}(a) and (b) show the OFM and COFM respectively for the unstretched caldera potential.
Here the classifications are restricted to the two possible reactive behaviours for trajectories crossing this PODS -- right to left dynamical matching, and right to right broken dynamical matching.
The OFM [Fig.~\ref{fig:PODS}(a)] clearly shows the two natural behaviours for the PODS: Reactive trajectories coming from channel 2 and exiting across the caldera via channel 4, and ``anti-reactive'' trajectories tracing out the opposite journey, starting from channel 4 and leaving through channel 2.
While not shown here, if the other half of the surface with $p_y>0$ is included, exactly half of this full surface corresponds to $2-4$ dynamics, and the other half to $4-2$.

Now stretching the caldera past the critical value, to $\lambda=0.75$, we note the appearance of other dynamics on the PODS in Figs.~\ref{fig:PODS}(c) and (d).
Thus we see that the intersection of manifolds seen on the constant $y$ surface of section in Fig.~\ref{fig:SC_lobes} directly implies the existence of manifold lobes on the PODS.
If we think slightly differently about its definition, the PODS allows us to consider the surface in phase space where one pair of stable and unstable manifolds intersect (of course corresponding to the particular UPO under study).
Consequently, the intersection of any other manifold with either of these two specific manifolds anywhere in the phase space will necessarily result in the ``invasion'' of the PODS by these lobes of different asymptotic behaviours.
However, by restricting ourselves to this PODS, we automatically exclude behaviours that do not relate to the chosen channel, either through entrance or exit.
Thus we see that the complex yellow lobe appearing in the COFM of Fig.~\ref{fig:SC_lobes}(d) corresponds to the yellow lobe in Fig.~\ref{fig:PODS}(d).

It is worth emphasising that the very simplicity of these OFM plots on the PODS in Fig.~\ref{fig:PODS} reflects the significant advantage of considering the dynamics on such a PODS.
By essentially choosing a ``good coordinate system'', we are able to see relevant behaviours of interest (and those behaviours only) very clearly.
A comparison of Fig.~\ref{fig:SC_full}(a) with Fig.~\ref{fig:PODS}(a), and Fig.~\ref{fig:SC_fullstretch} with Fig.~\ref{fig:PODS}(c), makes this point clearly -- while these figures show the same system at the same (or similar in the stretched case) parameters, the almost overwhelming complexity of the constant $y$ surface of section is completely eliminated in the PODS representation.
Note that while the $\lambda$ values are slightly different for Fig.~\ref{fig:SC_fullstretch} and Fig.~\ref{fig:PODS}(c), these show the same qualitative regimes.
Consequently we can see that to get a snapshot of the full dynamics of the system, we should consider some non-specific surface of section such as in Fig.~\ref{fig:SC_full}(a), while the PODS is much more useful as a blank canvas on which to illustrate the development of complex dynamics resulting from the variation of parameters.
In the system studied here, this difference is between having to carefully examine Fig.~\ref{fig:SC_fullstretch} for the lobe formation, and seeing the breaking of dynamical matching clearly appear on the PODS in Figs.~\ref{fig:PODS}(c) and (d).

\begin{figure}[tb]
    \centering
    \includegraphics[width=\columnwidth]{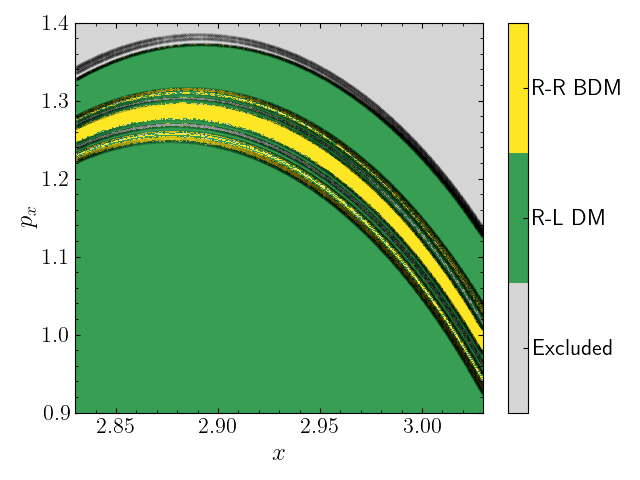}
    \caption{A zoomed in perspective of the PODS depicted in Fig.~\ref{fig:PODS}(d), associated with the top right saddle with $E=29$ and $\lambda=0.75$, showing the classification origin-fate map. The map is overlaid with the intersection of stable and unstable manifolds with the PODS, aligning exactly with the different trajectory behaviours.}
    \label{fig:SC_zoom}
\end{figure}

On the PODS, we are also able to directly compare the lobes shown by the OFM and those formed by the stable and unstable manifolds computed through the Lagrangian descriptors, in the same way as for the surface of section in Fig.~\ref{fig:SC_lobes}.
Figure~\ref{fig:SC_zoom} shows the COFM on the PODS associated with the TRS, overlaid with the manifolds in black.
For clarity, the figure zooms in on the relevant boundary region where the broken dynamical matching behaviour becomes evident.
Here the correspondence between the two methods for finding these transport lobes is clear once again, with the manifolds outlining the different regions identified by the OFM.
The fractal boundary of the lobes are visible in both the manifold and OFM visualisations, but the depth of the fractality is much more readily apparent from the OFM.
In addition, finding these higher-order manifolds responsible for the more complex transport mechanisms requires longer time Lagrangian descriptor computations (or extremely difficult eigenvector extrapolations from higher order periodic orbits), whereas a short time computation from the OFM is sufficient to completely demonstrate this phenomenon.
This confirms the utility of the OFM as both an investigative tool in its own right, as well as in a complementary role to manifold-based approaches.

\begin{figure}[tb]
    \centering
    \includegraphics[width=0.98\columnwidth]{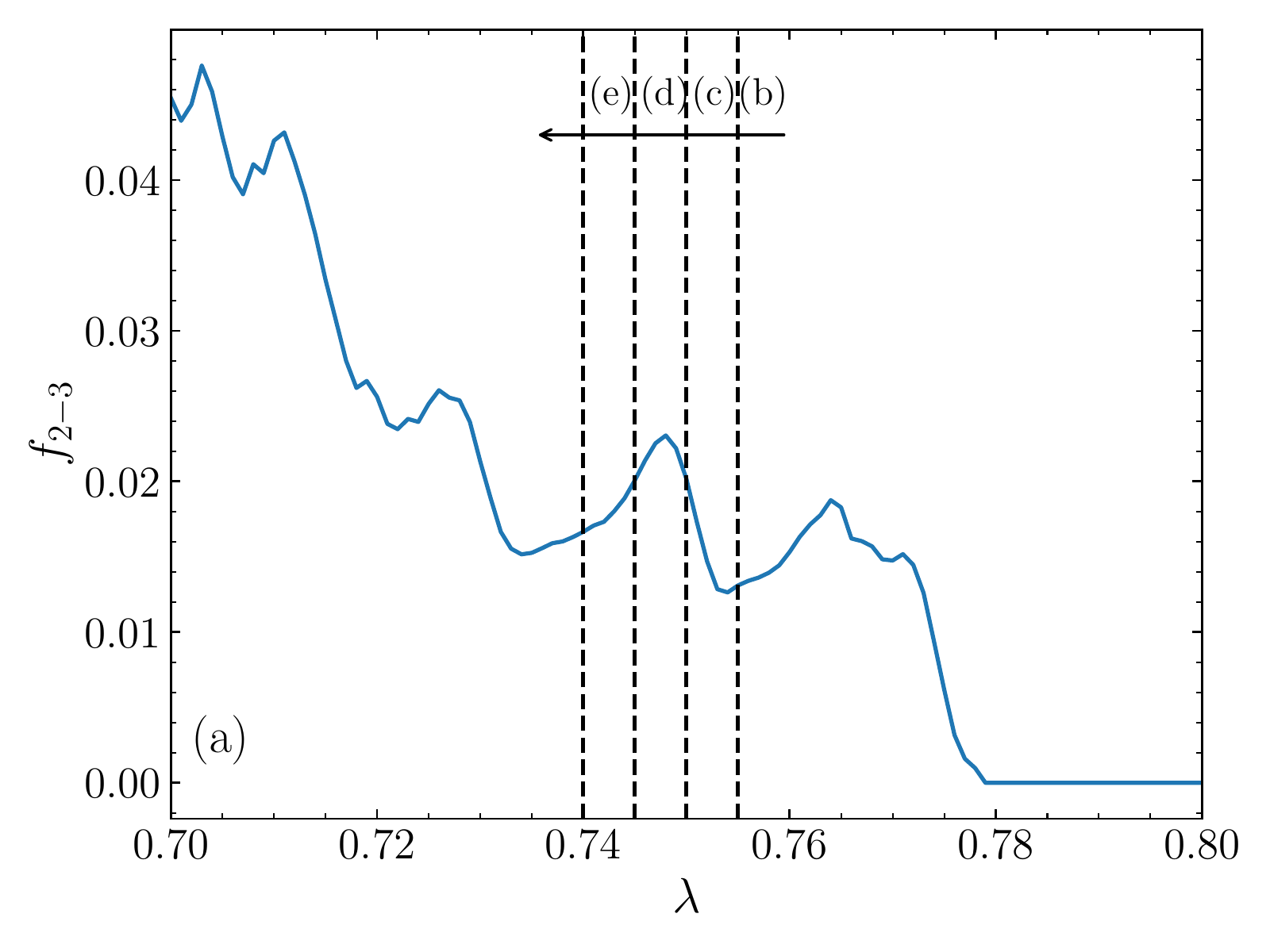}
    \includegraphics[width=0.45\columnwidth]{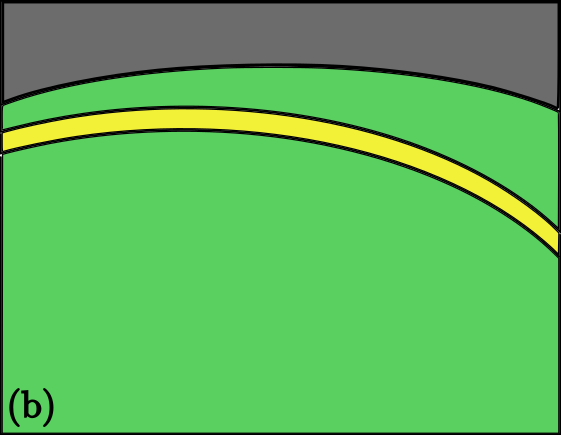}\hspace{0.2cm}
    \includegraphics[width=0.45\columnwidth]{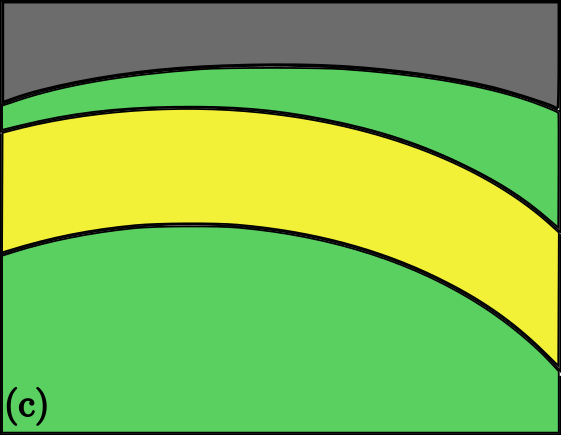}\vspace{0.2cm}
    \includegraphics[width=0.45\columnwidth]{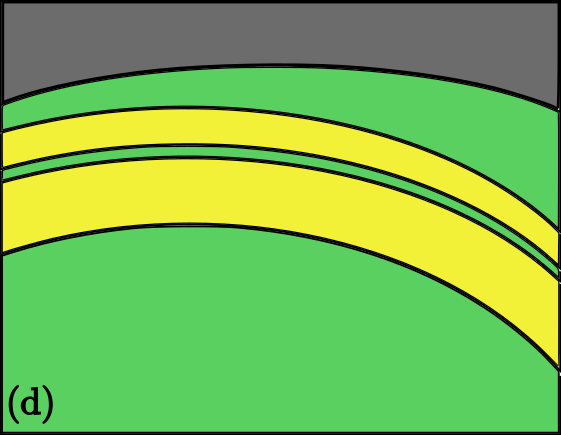}\hspace{0.2cm}
    \includegraphics[width=0.45\columnwidth]{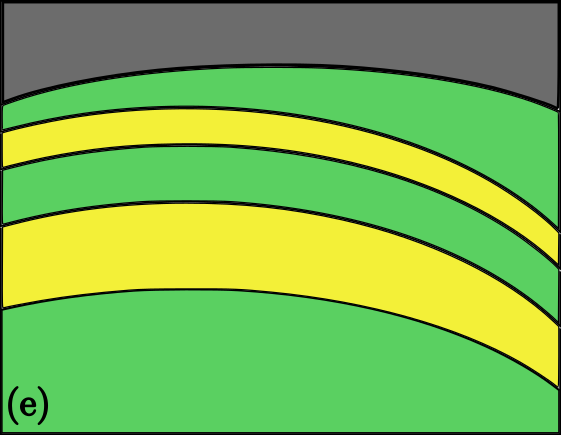}
    \caption{(a) The fraction $f_{2-3}$ of reactive trajectories crossing the PODS associated with the top right saddle that exit through the bottom right (breaking dynamical matching), as the stretching $\lambda$ varies. (b)-(e) A schematic representation of lobes appearing on the PODS that correspond to the increase and decrease in $f_{2-3}$ as $\lambda$ varies. Each diagram corresponds qualitatively to the points with labelled dashed lines in (a).}
    \label{fig:fractions}
\end{figure}

\subsection{Quantifying Origin-Fate Fractions} 
\label{sub:quantifying_origin_fate_fractions}

As a further application of the OFM on the PODS, here we demonstrate the quantitative prospects of this technique in the context of branching ratios.
The plots presented to this point amply show that we can detect the onset of broken dynamical matching by the appearance of particular lobes on a given surface of section.
However, as the OFM allocates an identifier to each initial condition studied corresponding to its dynamical behaviour, we can easily extract the distribution of origin-fate indices for the entire PODS, and particularly the branching ratio of reactants emanating from the top right.
In Fig.~\ref{fig:fractions}(a) we see the fraction of reactive (emanating from channel 2) trajectories corresponding to broken dynamical matching (i.e.~$2-3$) as the caldera is stretched by a factor of $\lambda$.
Viewed from right to left, as $\lambda$ is decreased from the neutral state, there is a clear transition from the dynamical matching regime of $f_{2-3}=0$ to the mixed regime with $f_{2-3}\ne0$ as the critical value of $\lambda\approx0.778$ is crossed.
This fraction is the quantitative representation of the area of the yellow lobe we see appearing in Fig.~\ref{fig:PODS}(d) (although we note that the fraction is computed using the entire PODS, not just this subsection of the $x-p_x$ projection), and the critical transition from zero to nonzero area corresponds to the intersection of the manifolds and the formation of the lobe on the OFM.
Note however that this use of the OFM allows us to detect this transition point with arbitrary accuracy, without requiring manifold computations.

Moving to this quantitative approach allows for a much more direct investigation of the effect of varying the stretching $\lambda$.
Particularly, an aspect which would not be clear from pure inspection of the OFM plots at different $\lambda$ values is the fluctuation of the area of the lobes, which appears in Fig.~\ref{fig:fractions}(a) as oscillations in the fraction of $2-3$ trajectories.
As we follow $f_{2-3}$ with decreasing $\lambda$, we see regions of increasing fraction of $2-3$ trajectories, as well as regions where $f_{2-3}$ decreases again.
We can then return to the OFM and look immediately at particular points of interest to understand the mechanism of these fluctuations.
Considering $\lambda$ values on either side of the peak near $\lambda=0.75$ [shown by dashed lines labelled (b)-(e) in Fig.~\ref{fig:fractions}(a)], we represent slightly simplified schematic versions of the COFM on the PODS at each of these stages.
We see that the increase in fraction of $2-3$ trajectories corresponds to a simple increase in the main area of the $2-3$ lobe, evidenced by the growth seen going from panel (b) to (c).
This growth does not continue smoothly however, but is interrupted by the appearance of another $2-4$ lobe inside the main $2-3$ lobe [Fig.~\ref{fig:fractions}(d)], the onset of which marks the beginning of the decrease in $f_{2-3}$.
This $2-4$ lobe grows in turn [see the transition from panels (d) to (e) in Fig.~\ref{fig:fractions}], reducing the fraction of $2-3$ trajectories until eventually another $2-3$ lobe appears and starts to grow itself.
This fractal cascade of lobes is repeated and corresponds to each of the phases of the increase and decrease of $f_{2-3}$.
Clearly the $2-3$ lobes are overall growing more quickly, as $f_{2-3}$ trends upwards as $\lambda$ decreases, but there is a continuous evolution of these very complex lobes on the PODS corresponding to the underlying nontrivial dynamics.
Apart from being dynamically a very interesting phenomenon, the physical interpretation of this fractal cascade would be that as the mechanism governing the caldera potential stretches it, the branching ratio will change unpredictably with small changes in the stretching.
A chemical reaction described by this potential would consequently yield different products for reactants with extremely similar initial states.

While this is but a single example, we emphasise that the OFM allowed the direct study of the phase space transport mechanisms without having to resort to an array of complex machinery.
It enables both the detection and the visualisation of this process, both of which would be virtually impossible via pure manifold considerations, or without the use of the PODS.
In particular, the flexibility of the OFM method enables us to simply quantify and explain a variety of phenomena at any parameter values of the system at hand.


\section{Summary and Conclusions} 
\label{sec:summary_and_conclusions}
The OFM is a simple and effective tool for classifying initial conditions in any kind of ``reactant-product'' system in terms of their long-time transport characteristics.
Extending the notion of basins of attraction, the OFM uses backward time integration to identify the origin of a given initial condition, and forward time integration to find its final fate.
Indexing each point on a given surface of section according to this origin-fate pairing enables a detailed understanding of lobe dynamics, particularly in conjunction with conventional manifold identification techniques.
We showed in Figs.~\ref{fig:SC_full}-\ref{fig:trajs} that the OFM is able not only to reproduce results obtained using such manifold techniques in a 2 degrees of freedom caldera potential, but provides a much more intricate description of the higher order dynamics with simpler computations.
An additional advantage of the method is that it allows for the accurate estimation of positions of unstable periodic orbits on a surface of section in conservative systems, due to the necessary distinct four-quadrant behaviours defined by the manifolds of this orbit (Fig.~\ref{fig:UPO_identify}).

As an application of the OFM, we investigated the phenomenon of dynamical matching, i.e.~having a one-to-one correspondence between origin and fate channels, in an symmetric caldera potential, with particular note given to the possibility of simply computing the OFM on a specific surface such as the PODS corresponding to a transition state.
By exploring the transport dynamics on a PODS associated with one entrance channel, we show that the critical stretching corresponding to the breaking of dynamical matching is associated with the formation of complex fractal lobes on the PODS created by intersections of manifolds (Figs.~\ref{fig:PODS} and \ref{fig:SC_zoom}).
Restricting the map to only specific behaviours (origin-fate pairs) of interest, we can plot a COFM which succinctly displays the relevant information, clearly laying out the transport dynamics.

Furthermore, the origin-fate indexing allows for direct quantification of branching ratios and similar measurements.
We demonstrated in Fig.~\ref{fig:fractions} that this quantification allows for studying the effect of various parameters on the system, and consequently the detailed analysis of relevant cases.
In particular, this allowed us to identify a complex cascade of fractal lobes as the mechanism producing fluctuations in the branching ratio of the stretched caldera potential.

We believe this method has a diverse range of possible applications, and could be useful in any system with a notion of origin and fate states, be they attractors in dissipative systems, potential wells, or exit channels in open conservative systems.
Computing the OFM on well-chosen targeted surfaces, such as the PODS, further enhances the power of the method for investigating particular behaviours.
\section*{Acknowledgements}
M.H.~acknowledges funding from the National Research Foundation (NRF) of South Africa, Grant No.~129630, and was supported by the Max Planck Society.
M.H~and Ch.S.~thank the Center for High Performance Computing of South Africa for providing their computational resources.
M.K.~and S.W.~acknowledge the financial support provided by the EPSRC Grant No. EP/P021123/1.

\bibliography{HKWS}
\end{document}